\documentclass[11pt]{article}
\usepackage[margin=1in]{geometry}

\usepackage[utf8]{inputenc}
\usepackage{lmodern}
\usepackage[normalem]{ulem} 
\usepackage{float}
\usepackage{enumitem}
\usepackage[utf8]{inputenc} 
\usepackage{tocloft}
\usepackage{graphicx, titlesec}
\graphicspath{/figures}
\usepackage{amscd,amsmath,amsthm,amssymb}
\usepackage{verbatim}
\usepackage[all]{xy}
\usepackage{color}
\usepackage{url}
\usepackage{tikz}
\usepackage{mathtools}
\usepackage{mdframed, xcolor}
\usepackage[]{hyperref}

\hypersetup{
    colorlinks=true,
    linkcolor=violet,
    filecolor=black,      
    urlcolor=violet,
    citecolor=violet,
    pdfpagemode=FullScreen,
    }

\newtheorem{theorem}{Theorem}

\newtheorem{lemma}[theorem]{Lemma}

\newtheorem{proposition}[theorem]{Proposition}

\theoremstyle{definition}
\newtheorem{remark}[theorem]{Remark}
\newtheorem{example}[theorem]{Example}

\newcommand{\diag}{\operatorname{diag}}
\newcommand{\rank}{\operatorname{rank}}

\newcommand{\im}{\operatorname{im}}

\newcommand{\Spec}{\operatorname{Spec}}

\newcommand{\M}{\mathcal{M}}

\newcommand{\RR}{\mathbb{R}}
\newcommand{\CC}{\mathbb{C}}
\newcommand{\PP}{\mathbb{P}}
\newcommand{\bbS}{\mathbb{S}}

\newcommand{\Sym}{\mathrm{Sym}}
\newcommand{\GL}{\mathrm{GL}}
\newcommand{\variety}[1]{\mathcal{V}(#1)}
\newcommand{\orbit}[1]{\mathcal{O}(#1)}
\newcommand{\cA}{\mathcal{A}}
\newcommand{\cB}{\mathcal{B}}
\newcommand{\cV}{\mathcal{V}}

\newcommand{\id}{\mathrm{id}}
\newcommand{\LT}{\mathrm{LT}}

\newif\ifincludesupplement
\includesupplementfalse 

\newcommand{\authorinfo}{%
  \par\vspace{2em}
  \noindent
  Yulia Alexandr\\
  School of Mathematics and Statistics,
  University of Melbourne\\
  \href{mailto:yulia.alexandr@unimelb.edu.au}
       {\texttt{yulia.alexandr@unimelb.edu.au}}

  \par\bigskip\noindent
  Moritz Grillo\\
  Max Planck Institute for Mathematics in the Sciences\\
  \href{mailto:Moritz.Grillo@mis.mpg.de}
       {\texttt{moritz.grillo@mis.mpg.de}}

  \par\bigskip\noindent
  Thiago Holleben\\
  Max Planck Institute for Mathematics in the Sciences\\
  \href{mailto:Thiago.Holleben@mis.mpg.de}
       {\texttt{thiago.holleben@mis.mpg.de}}

  \par\bigskip\noindent
  Guido Mont\'ufar\\
  Department of Mathematics, Department of Statistics \& Data Science,
  UCLA\\
  Max Planck Institute for Mathematics in the Sciences\\
  \href{mailto:montufar@math.ucla.edu}
       {\texttt{montufar@math.ucla.edu}}

  \par\bigskip\noindent
  Teresa Yu\\
  Department of Mathematics,
  University of Notre Dame\\
  \href{mailto:tyu6@nd.edu}
       {\texttt{tyu6@nd.edu}}
}

\title{Defining equations for two-block shallow ReLU pattern varieties}
\author{Yulia Alexandr, Moritz Grillo, Thiago Holleben, \\ Guido Mont\'ufar, Teresa Yu}
\date{}

\begin{document}

\maketitle

\vspace{-1em}
\begin{abstract}
We give a complete algebraic characterization of two-block shallow ReLU pattern varieties. 
We show that the determinantal equations arising from rank constraints imposed by the activation patterns generate the defining ideals, thereby proving a conjecture of Alexandr--Mont\'ufar. 
The set-theoretic proof uses a characterization of orbit closures of matrix pencils in terms of the Kronecker canonical form. 
\end{abstract}

\maketitle

\section{Introduction}\label{sec:intro}

\subsection{Overview}

A ReLU neural network defines a piecewise linear map from inputs to outputs. 
Within each input region where the same hidden neurons are active, the network acts as a linear map. Since the different linear maps share the same network parameters, their coefficients satisfy algebraic relations. 
Describing these relations helps characterize the functions represented by a given ReLU architecture. 
\emph{ReLU pattern varieties} were introduced in~\cite{alexandr_montufar_2027} to study these algebraic relations. 

For a fixed collection of activation patterns, the associated pattern variety is the Zariski closure of the image of the polynomial map that sends network parameters to the corresponding linear maps. 
For a pair of activation patterns in a shallow network, natural rank constraints give determinantal equations that vanish on the associated pattern variety. 
These equations were conjectured to generate the defining ideal \cite[Conjecture 26]{alexandr_montufar_2027}. 
In this paper, we prove this conjecture for arbitrary input, hidden, and output dimensions. 

To state the main result, we give an equivalent definition of the two-block pattern variety. 
Fix integers $n_0,n_1,n_2\geq 1$ and $a,b,c\geq 0$ with $a+b+c=n_1$. 
Here $a,b,c$ count the hidden neurons active only in the first pattern, only in the second pattern, and in both patterns, respectively. 
The associated ReLU pattern variety $\mathcal{U}_\cA \subseteq \CC^{n_2\times n_0}\times\CC^{n_2\times n_0}$ is the Zariski closure of the following set of matrix pairs: 
\begin{equation}\label{eqn:join defn}
    \{(A+C,B+C):\rank A\leq a,\,\rank B\leq b,\,\rank C\leq c\}.
\end{equation}

\begin{theorem}\label{thm:main}
    Let $R=\CC[X,Y]$ denote the coordinate ring of $\CC^{n_2\times n_0}\times\CC^{n_2\times n_0}$. Then the defining prime ideal of the ReLU pattern variety $\mathcal{U}_\cA$ is generated by the following minors:
    \begin{enumerate}
        \item $(a+c+1)$-minors of $X$;
        \item $(b+c+1)$-minors of $Y$;
        \item $(n_1+1)$-minors of $[X\mid Y]$ and $[X^\top\mid Y^\top]$;
    \item $(a+b+1)$-minors of $X-Y$.
    \end{enumerate}
\end{theorem}

For a shallow ReLU network without biases, the linear map associated with an activation pattern is the sum of the rank-one contributions of its active hidden neurons. 
If we consider the pattern variety associated to \emph{all} full-dimensional activation regions, its ideal is generated by linear equations and $2$-minors expressing these rank-one constraints; see~\cite{finkel2026canonicalactivationpatternsrelu} and Section~\ref{subsec:more act regions}
for further~discussion.

Selecting a subset of activation patterns amounts to taking the Zariski closure of a coordinate projection of the full pattern variety. 
The defining ideal of the pattern variety for a subset of activation patterns therefore amounts to considering an elimination ideal, for which finding explicit generators is generally a difficult nonlinear algebra problem.
Theorem~\ref{thm:main} describes the generators for the defining ideal for the variety associated to any pair of activation patterns. In particular, these impose polynomial constraints on the two observed maps associated to the pair of patterns,
whether or not the grouped neuron contributions can be recovered linearly from the two maps.
Our result also imposes necessary compatibility constraints on every pair of patterns' linear maps within a larger collection of patterns.

We first establish the set-theoretic version of the theorem by identifying our variety with the orbit closure of a specific matrix pencil. We show that a characterization of matrix pencil orbit closures in terms of the Weyr characteristics of their Kronecker invariants corresponds to the stated rank constraints. To prove the full ideal-theoretic statement, we show that these set-theoretic equations generate a prime ideal. 
We first establish primality in the square case where the input, hidden, and output layers have equal width, $n_0=n_1=n_2$, 
by proving Cohen--Macaulayness and generic reducedness. 
Restriction and inheritance arguments then extend primality to arbitrary input and output dimensions. 
We also give an alternative proof of both the set-theoretic and ideal-theoretic statements in the square case using an incidence variety over a product of Grassmannians. 

The image of the parametrization map defining the pattern variety is not always Zariski closed. 
We give examples of pairs of matrices that are limits of pairs in the image, but are not themselves attained by the parametrization. 
Proposition~\ref{prop:closed image} characterizes when the image is Zariski closed in terms of the architecture and activation patterns.

\subsection{Related work}

\subsubsection{ReLU networks}
There has been much interest in ReLU networks from algebraic, geometric, and combinatorial perspectives in recent years \cite{brandenburg2024the,flinth2026identifiability,grigsby2022transversality, masden2025algorithmic, montufar2014number, zhang2018tropical}, to name a few among many other works.

The aforementioned work~\cite{alexandr_montufar_2027} initiates the systematic study of algebraic relations among the linear pieces of ReLU networks through ReLU pattern varieties. 
These relations, which are the main focus of the present work, are naturally connected to functional dimension and parameter identifiability \cite{grigsby2025functional, stock2023embedding}. 
For a fixed finite dataset, the local dimension of the set of outputs obtained by varying the network parameters is studied in~\cite{bonapellissier2024geometry}.
Rank conditions and differences between linear pieces arise naturally in this literature. 
For example, rank-one differences between the linear parts on adjacent regions enter the geometric study of identifiability \cite{grillo2026most} and parameter fibers \cite{gegenfurtner2026symmetries}. 

Questions about the existence of best approximations by ReLU networks are closely related to whether the parametrization image is closed.  
For finite datasets, it was shown in~\cite{lim2022best} that the squared-error infimum need not be attained, even for shallow networks with vector-valued outputs. 
For population approximation of continuous scalar-valued functions, best approximations by shallow ReLU networks exist under suitable assumptions on the input measure~\cite{dereich2024}. 
For shallow scalar-valued ReLU networks of fixed width, the realization space is closed in the uniform norm on bounded cubes~\cite{petersen2021topological}.

\subsubsection{Secant varieties of Segre products}

Our defining equations connect ReLU pattern varieties to secant varieties of Segre varieties. 
Identifying a pair of $n_2\times n_0$ matrices $(X,Y)$ with the tensor $e_1\otimes X+e_2\otimes Y$, the decomposition $(X,Y)=(A+C,B+C)$ gives
\[
e_1\otimes X+e_2\otimes Y
=
e_1\otimes A+e_2\otimes B+(e_1+e_2)\otimes C.
\]
Decomposing $A,B,C$ into sums of rank-one matrices expresses this tensor as a sum of at most $a+b+c=n_1$ rank-one tensors. This places $\mathcal{U}_{\mathcal{A}}$ inside the affine cone over the $n_1$-th secant variety of
$\operatorname{Seg}
\bigl(
\mathbb{P}^{1}\times
\mathbb{P}^{n_2-1}\times
\mathbb{P}^{n_0-1}
\bigr)$. 
For this tensor format, the defining ideal of the secant variety is generated by the $(n_1+1)$-minors of the two concatenation flattenings~\cite[Theorem~1.1]{LW}. 
This description holds in arbitrary characteristic, where the flattening minors also form a Gr\"obner basis for suitable term orders~\cite[Theorem~4.1]{CDS}. 
The generators are precisely the concatenation minors in family~(3) of Theorem~\ref{thm:main}.
The remaining three families of minors impose rank bounds on $X$, $Y$, and $X-Y$. Theorem~\ref{thm:main} shows that adjoining these minors to the secant ideal gives the defining ideal of the pattern variety. 

\subsubsection{Matrix pencils}
Another connection, central to our set-theoretic proof, comes from viewing pairs of matrices as matrix pencils. 
The connection between tensors with a two-dimensional factor and the Kronecker canonical form is discussed in~\cite[Section~3.11.1]{landsberg}. 
Criteria for containment in orbit closures of matrix pencils were established in~\cite{pokrzywa1986} and subsequently formulated in terms of Weyr characteristics in~\cite{dehoyos1990}. 
We will use these criteria to derive the set-theoretic equations for the pattern variety.

\vspace{-0.5em}
\subsection{Structure of the paper}
The paper is organized as follows. Section~\ref{sec:prelim} introduces two-block ReLU pattern varieties and their natural determinantal equations, and characterizes when the parametrization image is closed. Section~\ref{sec:set-theoretic} establishes the set-theoretic description using matrix pencil orbit closures. Section~\ref{sec:ideal-theoretic} proves that the proposed minors generate the defining ideal, first in the square case and then in arbitrary dimensions. Section~\ref{sec:incidence} gives an alternative proof of primality in the square case using an incidence variety. We end in Section~\ref{sec:conclusion} with some directions for future work on defining ideals for more activation patterns and Gr\"obner bases for these ideals.

\vspace{-0.5em}
\subsubsection*{Acknowledgments}
We thank Ben Hollering and Mateusz Micha{\l}ek for helpful discussions. MG and GM were supported by the DFG under project 464109215 within the Priority Programme SPP 2298, ``Theoretical Foundations of Deep Learning''. 
YA and GM were supported by the DARPA Artificial Intelligence Quantified (AIQ) program under grant HR00112520014. 
GM was also supported in part by 
NSF grants 
2446222 (RTG), 
DMS-2522495 (AIMing), 
CCF-2212520 (MoDL), 
and DMS-2145630 (CAREER), 
and by the BMFTR through DAAD project 57616814 (SECAI).

\subsubsection*{AI Disclosure}
We used ChatGPT-5.6 to generate the main proof strategies for Sections~\ref{sec:matrix pencils} and~\ref{sec:ideal-theoretic}, including the connection to matrix pencil orbit closures. 
We then worked through the proposed strategies, verified the arguments, and worked out the mathematical details needed to develop these ideas into the proofs presented here.
We used ChatGPT-6.1 Astra to assist with the proofs of Lemma~\ref{lem:pi_2 partial G dim} and Proposition~\ref{prop:gb countereg}, as well as with revising, polishing, and proofreading the manuscript.

\section{Preliminaries}
\label{sec:prelim}\label{sec:preliminaries}

We introduce two-block ReLU pattern varieties and the natural rank constraints they satisfy. We also characterize when the parametrization image is Zariski closed.

\subsection{Definitions}

A \textit{shallow ReLU network} with architecture $(n_0,n_1,n_2)$ is a function
$$
f_\theta:\RR^{n_0}\to\RR^{n_2},
\qquad
f_\theta(x)=W^{(2)}\sigma(W^{(1)}x),
$$
where
$$
\theta = (W^{(1)}, W^{(2)}),\qquad  W^{(1)}\in\RR^{n_1\times n_0},
\qquad
W^{(2)}\in\RR^{n_2\times n_1},
$$
are parameters, and the \textit{activation function}
$\sigma(z)=\max\{0,z\}$ is applied coordinatewise. 

For fixed parameters $\theta$, the input space is partitioned into activation regions according to the \textit{activation pattern}
$$
A(x)=(A_1(x),\ldots,A_{n_1}(x))\in\{0,1\}^{n_1},
$$
where $A_i(x)=1$ if $(W^{(1)}x)_i>0$ and $A_i(x)=0$ otherwise.
For $A\in\{0,1\}^{n_1}$, the corresponding \textit{activation region} is
$$
R_A=\{x\in\RR^{n_0}:A(x)=A\}.
$$
On $R_A$, the network is a linear function of the data $x$: 
$$
f_\theta(x)=M_A(\theta)x,
\qquad x\in R_A, 
\qquad
M_A(\theta)=W^{(2)}\diag(A)W^{(1)}.
$$

We now fix two activation patterns $\mathcal A=(A_1,A_2)$ and write
$M_i(\theta):=M_{A_i}(\theta)$ for $i=1,2$. Varying the network
parameters gives the polynomial map:
$$
\varphi_{\mathcal A}:\theta\longmapsto
\big(M_1(\theta),M_2(\theta)\big).
$$
The same polynomial formulas define $\varphi_{\mathcal A}$
on the complex parameter space
$\CC^{n_1\times n_0}\times\CC^{n_2\times n_1}$.
Henceforth, $\varphi_{\mathcal A}$ denotes this complexified map.
The \textit{two-block ReLU pattern variety} $\mathcal U_{\mathcal A}:=\overline{\im \varphi_{\mathcal A}}$ is
the Zariski closure over $\CC$ of the image of $\varphi_{\mathcal A}$. 
This construction records the linear maps associated with two formal activation patterns. 
For a fixed collection of activation patterns, there is only a semialgebraic subset in the parameter space where the parameters actually induce these activation patterns. We neglect this in our setup and allow arbitrary parameters, thereby not imposing the inequalities ensuring that the selected activation regions have nonempty interior.

Let $a=\sum_j (A_1)_j(1-(A_2)_j)$ be the number of neurons active only in the first pattern, let $b=\sum_j (1-(A_1)_j)(A_2)_j$ be the number active only in the second pattern, and let $c=\sum_j (A_1)_j(A_2)_j$ be the number active in both patterns. Since neurons that are inactive in both patterns do not contribute to either $M_1$ or $M_2$, we may discard them and assume that every hidden neuron is active in at least one pattern. Hence $n_1=a+b+c$. 
After reordering the hidden neurons accordingly, we may write
$$
W^{(2)}=[\,U_1\ U_2\ U_{12}\,],
\qquad
W^{(1)}=
\begin{bmatrix}
V_1\\
V_2\\
V_{12}
\end{bmatrix},
$$
where the three blocks have sizes $a,b,c$. Then
$$
M_1=U_1V_1+U_{12}V_{12},
\qquad
M_2=U_2V_2+U_{12}V_{12}.
$$
Setting $A=U_1V_1$, $B=U_2V_2$, and $C=U_{12}V_{12}$, we get
$$M_1=A+C,\quad M_2=B+C\quad\text{with}\quad \rank A\leq a,\quad \rank B\leq b,\quad \rank C\leq c.$$
In particular,
\begin{equation}\label{eq:rank-conditions}
\rank M_1\leq a+c,\qquad
\rank M_2\leq b+c,\qquad
\rank(M_1-M_2)\leq a+b.
\end{equation}
Moreover, the column spaces of $M_1$ and $M_2$ are both contained in the
column space of $W^{(2)}$, while their row spaces are both contained in the
row space of $W^{(1)}$. Hence
\begin{equation}\label{eq:support-rank-conditions}
\rank[\,M_1\mid M_2\,]\leq n_1,
\qquad
\rank
\begin{bmatrix}
M_1\\
M_2
\end{bmatrix}
\leq n_1.
\end{equation}

Let $M_1=(x_{jk})$ and $M_2=(y_{jk})$ be generic $n_2\times n_0$ matrices, and let
$R=\CC[x_{jk},y_{jk}]$ be the coordinate ring of
$\CC^{n_2\times n_0}\times\CC^{n_2\times n_0}$. We also write $X=M_1$ and $Y=M_2$, as in the introduction. Let
$J_{\mathcal A}=I(\mathcal U_{\mathcal A})\subseteq R$ be the defining
prime ideal of the two-block ReLU pattern variety. For a matrix $N$, let $I_k(N)$
denote the ideal generated by its $k\times k$ minors; we set $I_k(N)=0$ if $k$ is greater than either of the dimensions of $N$. Define
\begin{equation} 
\begin{aligned}
I_{\mathcal A}
={}&I_{a+c+1}(M_1)+I_{b+c+1}(M_2)
   +I_{a+b+1}(M_1-M_2)
+I_{n_1+1}([\,M_1\mid M_2\,])
   +I_{n_1+1}\left(
   \begin{bmatrix}
   M_1\\
   M_2
   \end{bmatrix}\right).
\end{aligned}
\label{eq:I_A}
\end{equation} 
By \eqref{eq:rank-conditions} and \eqref{eq:support-rank-conditions}, every generator of $I_{\mathcal A}$ vanishes on $\mathcal U_{\mathcal A}$, and hence $I_{\mathcal A}\subseteq J_{\mathcal A}.$
It was conjectured in~\cite{alexandr_montufar_2027} that this inclusion is an equality. We prove this conjecture.

\begin{example}
Consider the architecture \((4,3,4)\) with activation patterns \(A_1=(1,1,0)\) and \(A_2=(0,1,1)\), so that \(a=b=c=1\). For generic \(4\times4\) matrices \(M_1,M_2\), we have
\[
\begin{aligned}
I_{\mathcal A}={}&I_3(M_1)+I_3(M_2)+I_3(M_1-M_2)+I_4([M_1\mid M_2])
+I_4\!\left(\begin{bmatrix}M_1\\M_2\end{bmatrix}\right).
\end{aligned}
\]
All four families appear in the defining ideal. A minimal generating set consists of \(16 + 16 + 16 = 48\) cubics and \(20 + 20 = 40\) quartics.
\end{example}

\subsection{Closedness of the image}

The following example shows that the Zariski closure can contain pairs not attained by the parametrization.  
Below, we characterize precisely when the image is closed, or equivalently, $\im\varphi_{\mathcal A}=\mathcal U_{\mathcal A}$. 

\begin{example}
\label{ex:nonclosed-image}
Consider the architecture $(2,3,2)$ with $a=b=c=1$.
For $t\neq0$, let
\[
A_t=\begin{pmatrix}0&-t^{-1}\\0&1\end{pmatrix},
\qquad
B_t=\begin{pmatrix}1&-t^{-1}\\-2t&2\end{pmatrix},
\qquad
C_t=\begin{pmatrix}1&t^{-1}\\0&0\end{pmatrix}.
\]
Each of these matrices has rank one. Hence
\[
(M_{1,t},M_{2,t})
=(A_t+C_t,B_t+C_t)
=
\left(
I_2,\begin{pmatrix}2&0\\-2t&2\end{pmatrix}
\right)
\]
belongs to $\im\varphi_{\mathcal A}$ for $t\neq0$. 
Although the matrices $A_t,B_t$ and $C_t$ diverge as $t\to0$,
their sums converge:
$
(M_{1,t},M_{2,t})\longrightarrow(I_2,2I_2).$

This limit does not belong to the image of $\varphi_{\mathcal A}$. 
Otherwise, writing $(M_1,M_2)=(I_2,2I_2)$, there would exist
$U\in\CC^{2\times3}$ and $V\in\CC^{3\times2}$ such that
$
M_1=U\diag(1,0,1)V$ and $
M_2=U\diag(0,1,1)V.
$
Since $\rank M_1=2$, both $U$ and $V$ have rank two. 
However,
$
2M_1-M_2 = U\diag(2,-1,1)V,
$
and the middle matrix is invertible. 
On the other hand $2M_1-M_2=2I_2-2I_2=0$. 
Sylvester's rank inequality therefore gives
$
0=\rank(2M_1-M_2)
\geq\rank U+\rank V-3
=1,
$
a contradiction.
Thus $\im\varphi_{\mathcal A}$ is neither Euclidean closed
nor Zariski closed.
\end{example}

\begin{proposition}\label{prop:closed image}
Let $n_1= a+b+c$ and $r= \min \{n_0,n_2\}$. Then, 
    $\im \varphi_\mathcal{A}$ is closed if and only if at least one of the following conditions holds:
    \begin{enumerate}
        \item  $abc =0$;
        \item $\max\{a,b,c\}\geq r$;
        \item $n_1\geq 2r$.
    \end{enumerate}
\end{proposition}

\begin{proof}
 We first prove closedness under each of the three conditions. Assume that $c=0$. Then, $\im \varphi_\mathcal{A}$ is given by the rank constraints $\rank M_1 \le a$ and $\rank M_2 \le b$. If $a=0$, then the image is given by the rank constraints $\rank M_1 \le c$ and $\rank (M_1- M_2) \le b$ and similarly if $b=0$ and hence $\im \varphi_\mathcal{A}$ is closed if $abc=0$.
    
    Next, if $a\ge r$, then $M_1$ is unconstrained and consequently $\im \varphi_\mathcal{A}$ is given by $\rank M_2 \leq b+c$. 
    Analogously, if $b\ge r$ (respectively $c \ge r$), the image $\im \varphi_\mathcal{A}$ is given by the rank constraints $\rank M_1 \leq a+c$ (respectively $\rank M_1-M_2 \leq a+b$).
    
    Now assume that $n_1 \geq 2r$ and by the preceding cases, we may assume $0<a,b,c<r$. We show that in that case, the image $\im \varphi_\mathcal{A}$ equals the whole ambient space. Transposing the matrices if necessary, assume that $n_2=r$. Set $p=r-a$ and $s=r-b$, then $p+s=2r-a-b\leq c<r$.
    For an arbitrary pair $(M_1,M_2)$, partition the rows into blocks of sizes $p,s,r-p-s$, and write
\[
M_1=\begin{bmatrix}X_1\\X_2\\X_3\end{bmatrix},
\qquad
M_2=\begin{bmatrix}Y_1\\Y_2\\Y_3\end{bmatrix}.
\]
Choosing the matrices
\[
C=\begin{bmatrix}X_1\\Y_2\\0\end{bmatrix},
\qquad
A=M_1-C,\qquad B=M_2-C
\] 
yields 
\[
\rank A\leq r-p=a,\qquad
\rank B\leq r-s=b,\qquad
\rank C\leq p+s\leq c,
\]
and thus $(M_1,M_2)\in\im\varphi_{\mathcal A}$.

For the opposite direction, assume that 1--3 do not hold.  
We may assume that $a=\max\{a,b,c\}$, since swapping $M_1,M_2$ interchanges $a,b$, while replacing
$(M_1,M_2)$ by $(M_1,M_1-M_2)$ interchanges $a,c$ and these invertible linear transformations preserve closedness. Set $d=\min\{r,a+b,a+c\}$.
Since $a=\max\{a,b,c\}$ and $b,c>0$, both $a+b$ and $a+c$
are greater than $a$ and $n_1/2$.
Together with our assumptions on $r$, this gives $d>a$
and $2d>n_1$.
Define $u=a+b-d$, $v=a+c-d$, $w=d-a-1$, and $s=2d-n_1-1$.
These are nonnegative integers satisfying $2+u+v+w+s=d$.
Let $A_t^0,B_t^0,C_t^0$ denote the $2\times2$ matrices from
Example~\ref{ex:nonclosed-image}. For $t\neq0$, form
\begin{align*}
A_t=A_t^0\oplus I_u\oplus(-I_v)\oplus0_w\oplus I_s,\quad
B_t=B_t^0\oplus2I_u\oplus0_v\oplus I_w\oplus0_s,\quad
C_t=C_t^0\oplus0_u\oplus2I_v\oplus I_w\oplus0_s.
\end{align*}
  These matrices have ranks
\begin{align*}
\rank A_t&=1+u+v+s=a,\quad
\rank B_t=1+u+w=b,\quad
\rank C_t=1+v+w=c.
\end{align*}
Moreover, $A_t+C_t=I_d$ and
\begin{equation*}
B_t+C_t=
\begin{pmatrix}2&0\\-2t&2\end{pmatrix}
\oplus2I_{u+v+w}\oplus0_s.
\end{equation*}
Since $d\leq r$, we can pad these matrices with zeros to size
$n_2\times n_0$. The resulting pairs  $(A_t+C_t,B_t+C_t)$ belong to
$\im\varphi_{\mathcal A}$ and  converge for $t \to 0$ to the zero-padded pair
$(M_1,M_2)=(I_d,2I_{d-s}\oplus0_s)$.
In particular, $\rank M_1=d$ and $\rank(2M_1-M_2)=s$.

Suppose that the limit $(M_1,M_2)$ belonged to the image. Then there would exist $U\in\CC^{n_2\times n_1}$ and
$V\in\CC^{n_1\times n_0}$ such that
$
M_1=U\diag(I_a,0_b,I_c)V $ and $
M_2=U\diag(0_a,I_b,I_c)V.$ Since $\rank M_1=d$, both $U$ and $V$ have rank at least $d$. The matrix $\diag(2I_a,-I_b,I_c)$ is invertible, so
Sylvester's rank inequality gives $ s=\rank(2M_1-M_2) = \rank\bigl(U\diag(2I_a,-I_b,I_c)V\bigr) \geq 2d-n_1=s+1$ yielding a contradiction. Thus, 
$\im\varphi_{\mathcal A}$ is not closed.
\end{proof}

When the image is not closed, the defining equations do not
distinguish attained pairs from those occurring only as limits.
Over complex parameters, the Euclidean and Zariski closures
of the image coincide, so every point of the pattern variety
is a limit of attained pairs. Exact representability requires a decomposition
$(X,Y)=(A+C,B+C)$ satisfying the prescribed rank bounds.
The existence of such a decomposition can be tested by
finitely many kernel computations and generic-rank tests;
we leave the development of this criterion to future work.

\section{Set-theoretic equations via matrix pencils}\label{sec:set-theoretic}
\label{sec:matrix pencils}

In this section, we prove the set-theoretic version of Theorem~\ref{thm:main}, i.e., we show that the proposed determinants (the generators of the ideal $I_{\mathcal A}$ described in \eqref{eq:I_A}) cut out the two-block ReLU pattern variety: $\mathcal U_{\mathcal{A}}=\variety{I_{\mathcal A}}$. 
Our proof translates this equality into a statement about orbit closures of matrix pencils. 
We first treat the case $n_1\leq\min\{n_0,n_2\}$, where the connection to matrix pencil orbits is most direct, and then extend the result to arbitrary dimensions.

\subsection{Matrix pencils and orbit closures}

For matrices $M_1,M_2\in\CC^{m\times n}$, we define the corresponding \textit{matrix pencil} to be
$$
M=xM_1+yM_2\in\CC[x,y]^{m\times n}.
$$
Thus $M$ is a matrix whose entries are homogeneous linear forms in $x$ and $y$. For
$\lambda=[\alpha:\beta]\in\PP^1$, we write
$$
M(\lambda)=\alpha M_1+\beta M_2.
$$
The \textit{normal rank} of $M$ is $\rank M=\rank_{\CC(x,y)}M.$ Equivalently, it is the rank of $M(\lambda)$ for a generic choice of $\lambda\in\PP^1$.
A point $\mu\in\PP^1$ is an \textit{eigenvalue} of $M$ if
$\rank M(\mu)<\rank M.$ Two matrix pencils $M$ and $L$ of the same dimensions are \textit{strictly equivalent} if $L=PMQ$
for some $P\in\GL_m$~and~$Q\in\GL_n$. The strict-equivalence orbit of $M$ is
$$
\orbit{M}=\{PMQ:P\in\GL_m,\ Q\in\GL_n\}.
$$
The Kronecker canonical form classifies matrix pencils up to strict equivalence~\cite[Section 3.11.1]{landsberg} (see also~\cite[Chapter XII, Sections 4--5]{gantmacher1959}). 
Every matrix pencil $M$ is strictly equivalent to a direct sum, unique up to reordering, of Jordan blocks, right singular blocks, and left singular blocks. For $\mu=[\alpha:\beta]\in\PP^1$, let
$\ell_\mu(x,y)=\beta x-\alpha y$, and choose a linear form
$q_\mu(x,y)$ such that $q_\mu(\mu)\neq 0$. A \textit{Jordan block} of size
$k$ associated with $\mu$ is the $k\times k$ pencil
$$
J_k(\mu)=\ell_\mu(x,y)I_k+q_\mu(x,y)N_k,
$$
where $N_k$ is the nilpotent Jordan block with ones on the superdiagonal.

Since $\det J_k(\mu)=\ell_\mu(x,y)^k$ and
$\ell_\mu(x,y)$ vanishes precisely when $[x:y]=\mu$, the evaluated block has rank $k$ whenever $[x:y]\neq\mu$. At $[x:y]=\mu$, it is a nonzero scalar multiple of $N_k$, because $q_\mu(\mu)\neq0$, and therefore has rank $k-1$.
For $k\geq 1$, a \textit{right singular block} is the $k\times(k+1)$ pencil
$$
R_k=
\begin{pmatrix}
x & y &        &        \\
  & x & y      &        \\
  &   & \ddots & \ddots \\
  &   &        & x & y
\end{pmatrix},
$$
and a \textit{left singular block} is its transpose $R_k^\top$. We also let
$R_0$ denote a zero column and $R_0^\top$ a zero row. The integer $k$ is
called the corresponding right or left \textit{minimal index}. It is the smallest degree of a nonzero polynomial relation
among the columns of $R_k$ or the rows of $R_k^\top$, respectively. 
Indeed,
$R_k$ has the polynomial kernel vector
$$
\big((-y)^k,x(-y)^{k-1},\ldots,x^{k-1}(-y),x^k\big)^\top,
$$
of degree $k$, and $R_k^\top$ has the same polynomial relation
on the left. Thus, in the Kronecker~canonical form of a pencil $M$, the Jordan blocks describe its rank drops at particular points of $\PP^1$, while the right and left singular blocks describe polynomial vectors in the kernels of $M$ and~$M^\top$, respectively. 

\begin{remark}
    In~\cite[Section 3.11.1]{landsberg} the Kronecker canonical form of a pencil $xM_1 + yM_2$ is described~as 
    $$
        xM_1 + yM_2 = \begin{pmatrix}
            R_{\varepsilon_1} \\
            & \ddots & \\
            & & R_{\varepsilon_q} \\
            & & &R_{\eta_1}^\top \\
            & & & &\ddots \\
            & & &&&R_{\eta_p}^\top \\
            &&&&&&xI_k + yF
        \end{pmatrix},
    $$
    where $F$ is a $k \times k$ matrix in Jordan canonical form. In this setting, the points $\lambda =[\alpha:\beta] \in \PP^1$ where the rank of $M(\lambda)$ drops are exactly the points where $\alpha = -\beta c$ for an eigenvalue $c$ of $F$. Note in particular that in this setting, the point $[1:0]$ is never an eigenvalue of a pencil. Since we can always apply a linear change of coordinates to move an eigenvalue away from $[1:0]$, we may view the description using the blocks $J_k(\mu)$ and the description with regular part $xI_k+yF$ as equivalent, up to applying linear changes of coordinates (see for example~\cite[Theorem 5.2]{BHMT2018}). 
\end{remark}

For an eigenvalue $\mu\in\PP^1$ of the pencil $M$, let $W(\mu,M)=(W_1(\mu,M),W_2(\mu,M),\ldots),$ where $W_i(\mu,M)$ is the number of Jordan blocks associated with $\mu$
whose size is at least $i$. Similarly, we define sequences
$r(M)=(r_0(M),r_1(M),\ldots)$ and $\ell(M)=(\ell_0(M),\ell_1(M),\ldots),$
where $r_i(M)$ and $\ell_i(M)$ are the numbers of right and left minimal
indices, respectively, that are at least $i$. These sequences are the
\textit{Weyr characteristics} of the corresponding Kronecker invariants. For two sequences $s=(s_1,s_2,\ldots)$ and $t=(t_1,t_2,\ldots)$ of nonnegative integers, we say that $t$ \textit{weakly majorizes}~$s$, and write $s\prec_W t$ if $s_1+\cdots+s_k\leq t_1+\cdots+t_k$ for every $k\geq 1$.

The following orbit closure result was first proved in~\cite[Theorem~3]{pokrzywa1986} and later reformulated in terms of Weyr characteristics in~\cite[Lemma~1.3]{dehoyos1990}; we use the formulation in~\cite[Theorem~1]{deteran2025dimensionorbitsmatrixpencils}. We write $\overline{\orbit{L}}^{\,E}$ for the Euclidean closure of $\orbit{L}$ and reserve $\overline{\orbit{L}}$ for its Zariski closure.

\begin{theorem}
\label{thm:orbit-closure-characterization}
Let $M$ and $L$ be matrix pencils of the same dimensions such that
$\rank M\leq\rank L$, and let
$$
\mathbf h=(h,h,\ldots), 
\qquad 
h=\rank L-\rank M 
.
$$
Then $M\in\overline{\orbit{L}}^{\,E}$ if and only if
\begin{enumerate}
\item $r(M)\prec_W r(L)+\mathbf h$;
\item $\ell(M)\prec_W \ell(L)+\mathbf h$;
\item $W(\lambda,L)\prec_W W(\lambda,M)+\mathbf h$ for every $\lambda\in\PP^1$.
\end{enumerate}
\end{theorem}

Since $\rank M\leq\rank L$, we have $h=\rank L-\rank M\geq 0$, so every entry of $\mathbf h$ is nonnegative. If $\lambda$ is not an eigenvalue of $L$, then $W(\lambda,L)=(0,0,\ldots)$, and condition~(3) holds automatically. Thus it is enough to check condition~(3) at the eigenvalues of $L$.

\begin{remark}
Although Theorem~\ref{thm:orbit-closure-characterization} is stated for Euclidean closure, the Euclidean and Zariski closures of $\orbit{L}$ coincide. Since every Zariski closed set is Euclidean closed, we have $\overline{\orbit{L}}^{\,E}\subseteq\overline{\orbit{L}}$.
Conversely, $\orbit{L}$ is Zariski open in its Zariski closure~\cite[Proposition in Section~1.8]{borelbook}. Since $\orbit{L}$ is Zariski dense in $\overline{\orbit{L}}$, its complement is a proper Zariski closed subset and hence has empty Euclidean interior. Thus $\orbit{L}$ is Euclidean dense in $\overline{\orbit{L}}$, and therefore
$\overline{\orbit{L}}^{\,E}=\overline{\orbit{L}}$.
\end{remark}

\subsection{The set-theoretic equality}

Assume that $n_1\leq\min\{n_0,n_2\}$. Define
$$
D_1=\diag(I_a,0_b,I_c),
\qquad
D_2=\diag(0_a,I_b,I_c),
\qquad 
D=xD_1+yD_2
=\diag(xI_a,yI_b,(x+y)I_c).
$$
Let
$$
L=
\begin{bmatrix}
D & 0\\
0 & 0
\end{bmatrix}
$$
be the $n_2\times n_0$ matrix pencil obtained from $D$ by padding with zero rows and columns. Let $\Psi: \CC^{n_2 \times n_0} \times \CC^{n_2 \times n_0} \to \CC[x,y]_1^{n_2 \times n_0}$ be the linear isomorphism defined by $\Psi(M_1,M_2)=xM_1+yM_2$, and set $X_{\mathcal A}
=
\Psi\bigl(\variety{I_{\mathcal A}}\bigr)$.

\begin{lemma}
\label{lem:equiv-statement-pencil}
Let $M_1,M_2\in\CC^{n_2\times n_0}$ and write
$\Psi(M_1,M_2)=xM_1+yM_2$. Let $L$ be the pencil defined above, and let
$X_{\mathcal A}=\Psi(\variety{I_{\mathcal A}})$.
Under the assumption $n_1\leq\min\{n_0,n_2\}$,
$
\mathcal U_{\mathcal A}=\variety{I_{\mathcal A}}
$
if and only if
$
\overline{\orbit{L}}=X_{\mathcal A}.
$
\end{lemma}

\begin{proof}
With the ordering of the hidden neurons fixed in Section~\ref{sec:preliminaries}, the first $a$ neurons are active only in the first pattern, the next $b$ only in the second, and the final $c$ in both. Therefore
$\diag(A_1)=D_1$ and $\diag(A_2)=D_2$. Hence
$$
\Psi(\im\varphi_{\mathcal A})
=
\{W^{(2)}DW^{(1)}:
W^{(2)}\in\CC^{n_2\times n_1},\
W^{(1)}\in\CC^{n_1\times n_0}\} = \{PLQ:
P\in\CC^{n_2\times n_2},\
Q\in\CC^{n_0\times n_0}\}.
$$

Indeed, $W^{(2)}$ can be chosen as the first $n_1$ columns of $P$, and
$W^{(1)}$ as the first $n_1$ rows of $Q$, and conversely every choice of
$W^{(2)}$ and $W^{(1)}$ can be extended in this way to matrices $P$ and $Q$.

Since $\GL_{n_2}\times\GL_{n_0}$ is dense in
$\CC^{n_2\times n_2}\times\CC^{n_0\times n_0}$, it follows that
$
\overline{\Psi(\im\varphi_{\mathcal A})}
=
\overline{\orbit{L}}.
$
Since $\Psi$ is a linear isomorphism, it preserves Zariski closures, so $\Psi(\mathcal U_{\mathcal A})=\Psi(\overline{\im\varphi_{\mathcal A}})=\overline{\Psi(\im\varphi_{\mathcal A})}=\overline{\orbit{L}}.$
On the other hand, by definition,
$\Psi(\variety{I_{\mathcal A}})=X_{\mathcal A}$. Hence
$\mathcal U_{\mathcal A}=\variety{I_{\mathcal A}}$ if and only if
$\overline{\orbit{L}}=X_{\mathcal A}$.
\end{proof}

We already know from Section~\ref{sec:preliminaries} that
$\mathcal U_{\mathcal A}\subseteq\variety{I_{\mathcal A}}$. Applying the linear isomorphism $\Psi$ to this inclusion gives
$\Psi(\mathcal U_{\mathcal A})\subseteq\Psi(\variety{I_{\mathcal A}})$.
In the proof of Lemma~\ref{lem:equiv-statement-pencil}, we showed that
$\Psi(\mathcal U_{\mathcal A})=\overline{\orbit{L}}$, while by definition
$\Psi(\variety{I_{\mathcal A}})=X_{\mathcal A}$. 
Therefore $\overline{\orbit{L}}\subseteq X_{\mathcal A}$.
Thus, to prove $\overline{\orbit{L}}=X_{\mathcal A}$, it remains only to show the reverse inclusion $X_{\mathcal A}\subseteq\overline{\orbit{L}}$.

\begin{lemma}\label{lem:rank-majorization}
Let $M=xM_1+yM_2\in X_{\mathcal A}$ and let
$h=\rank L-\rank M$. Then the first two conditions in
Theorem~\ref{thm:orbit-closure-characterization} are equivalent to the
support-rank conditions:
\begin{align*}
r(M)\prec_W r(L)+\mathbf h
&\iff
\rank
\begin{bmatrix}
M_1\\
M_2
\end{bmatrix}
\leq n_1,\\
\ell(M)\prec_W\ell(L)+\mathbf h
&\iff
\rank[\,M_1\mid M_2\,]\leq n_1.
\end{align*}
For the third condition, it suffices to check the three eigenvalues of $L$, where we get: 
\begin{align*}
W([1:0],L)\prec_W W([1:0],M)+\mathbf h
&\iff \rank M_1\leq a+c,\\
W([0:1],L)\prec_W W([0:1],M)+\mathbf h
&\iff \rank M_2\leq b+c,\\
W([1:-1],L)\prec_W W([1:-1],M)+\mathbf h
&\iff \rank(M_1-M_2)\leq a+b.
\end{align*}
\end{lemma}
\begin{proof}
Let $\rho=\rank M$. Since $\rank L=n_1$, we have
$h=n_1-\rho$. We first compute the Kronecker invariants of $L$. Since $L$ is obtained
from the generically invertible pencil
$$
D=\diag(xI_a,yI_b,(x+y)I_c)
$$
by appending $n_0-n_1$ zero columns and $n_2-n_1$ zero rows, all right
and left minimal indices of $L$ are zero. Hence
$$
r(L)=(n_0-n_1,0,0,\ldots),
\qquad
\ell(L)=(n_2-n_1,0,0,\ldots).
$$
Moreover, the eigenvalues of $L$ are $[1:0]$, $[0:1]$, and $[1:-1]$.
Since $L$ is diagonal, all corresponding Jordan blocks have size one,
and therefore
$$
W([1:0],L)=(b,0,0,\ldots),
\qquad
W([0:1],L)=(a,0,0,\ldots),
\qquad
W([1:-1],L)=(c,0,0,\ldots).
$$
For every other $\mu\in\PP^1$, we have $W(\mu,L)=(0,0,\ldots)$, so the
third condition in Theorem~\ref{thm:orbit-closure-characterization}
holds automatically. We next consider the right minimal indices of $M$. We have
$r_0(M)=n_0-\rho$, while $r_0(M)-r_1(M)$ is the number of right minimal
indices equal to zero. These correspond exactly to constant vectors
$v\in\CC^{n_0}$ in the polynomial kernel of $M$, that is,
$
M_1v=M_2v=0.
$
Thus
$$
\ker M_1\cap\ker M_2
=
\ker
\begin{bmatrix}
M_1\\
M_2
\end{bmatrix},
$$
and hence
$$
\rank
\begin{bmatrix}
M_1\\
M_2
\end{bmatrix}
=
n_0-(r_0(M)-r_1(M))
=
\rho+r_1(M).
$$
On the other hand,
$$
r(L)+\mathbf h
=
(n_0-\rho,h,h,\ldots).
$$
The first entries of $r(M)$ and $r(L)+\mathbf h$ agree. Since $r(M)$
is nonincreasing, the weak majorization condition is therefore equivalent
to $r_1(M)\leq h$. Consequently,
$$
r(M)\prec_W r(L)+\mathbf h
\iff
r_1(M)\leq n_1-\rho
\iff
\rank
\begin{bmatrix}
M_1\\
M_2
\end{bmatrix}
\leq n_1.
$$

The argument for the left minimal indices is analogous. We have
$\ell_0(M)=n_2-\rho$, and the left minimal indices equal to zero
correspond to
$
\ker M_1^\top\cap\ker M_2^\top.
$
Therefore
$
\rank[\,M_1\mid M_2\,]
=
\rho+\ell_1(M).
$
Since
$
\ell(L)+\mathbf h
=
(n_2-\rho,h,h,\ldots),
$
we obtain
$$
\ell(M)\prec_W\ell(L)+\mathbf h
\iff
\ell_1(M)\leq n_1-\rho
\iff
\rank[\,M_1\mid M_2\,]\leq n_1.
$$

It remains to consider the three eigenvalues of $L$. For any
$\mu\in\PP^1$, the first Weyr number of $M$ at $\mu$ is the rank drop
at $\mu$:
$$
W_1(\mu,M)=\rho-\rank M(\mu).
$$
Since each $W(\mu,L)$ has only its first entry nonzero, the
corresponding weak majorization condition is equivalent to the inequality
for the first partial sum. At $\mu=[1:0]$, we have $M(\mu)=M_1$, hence
\begin{align*}
W([1:0],L)\prec_W W([1:0],M)+\mathbf h
&\iff b\leq \rho-\rank M_1+h\\
&\iff \rank M_1\leq n_1-b\\
&\iff \rank M_1\leq a+c.
\end{align*}
Similarly,
\begin{align*}
W([0:1],L)\prec_W W([0:1],M)+\mathbf h
&\iff a\leq \rho-\rank M_2+h\\
&\iff \rank M_2\leq n_1-a\\
&\iff \rank M_2\leq b+c,
\end{align*}
and
\begin{align*}
W([1:-1],L)\prec_W W([1:-1],M)+\mathbf h
&\iff c\leq \rho-\rank(M_1-M_2)+h\\
&\iff \rank(M_1-M_2)\leq n_1-c\\
&\iff \rank(M_1-M_2)\leq a+b.
\end{align*}
\end{proof}

\begin{theorem}\label{thm:set-theoretic-small}
Suppose that $n_1\leq\min\{n_0,n_2\}$. Then
$
\mathcal U_{\mathcal A}=\variety{I_{\mathcal A}}.
$
\end{theorem}

\begin{proof}
By Lemma~\ref{lem:equiv-statement-pencil}, it suffices to show that
$\overline{\orbit{L}}=X_{\mathcal A}$. We already have
$\overline{\orbit{L}}\subseteq X_{\mathcal A}$. Conversely, let $M\in X_{\mathcal A}$. Then $\rank M\leq n_1=\rank L$,
and Lemma~\ref{lem:rank-majorization} shows that the rank conditions
defining $X_{\mathcal A}$ are exactly the nontrivial majorization
conditions in Theorem~\ref{thm:orbit-closure-characterization}.
Therefore $M\in\overline{\orbit{L}}$, proving the reverse inclusion.
\end{proof}

The assumption $n_1\leq\min\{n_0,n_2\}$ can be removed by padding the matrices $M_1$ and $M_2$ with zero rows and columns.

\begin{theorem}\label{thm:set-theoretic-general}
For arbitrary $n_0,n_1,n_2$, we have
$
\mathcal U_{\mathcal A}=\variety{I_{\mathcal A}}.
$
\end{theorem}

\begin{proof}
We already know that $\mathcal U_{\mathcal A}\subseteq\variety{I_{\mathcal A}}$, so it remains to prove the reverse inclusion. Let $(M_1,M_2)\in\variety{I_{\mathcal A}}$ and set $N=\max\{n_0,n_1,n_2\}$. Consider the enlarged architecture $(N,n_1,N)$ with the same pair of activation patterns $\mathcal A=(A_1,A_2)$. To distinguish the objects associated with this enlarged architecture, denote its parametrization, pattern variety, and determinantal ideal by $\varphi'_{\mathcal A}$, $\mathcal U'_{\mathcal A}$, and $I'_{\mathcal A}$, respectively. For $i=1,2$, let
$$
M_i'=
\begin{bmatrix}
M_i & 0\\
0 & 0
\end{bmatrix}
\in\CC^{N\times N},
$$
where the zero blocks have the appropriate sizes. Padding by zero rows and columns preserves the ranks of $M_1$, $M_2$, $M_1-M_2$, $[\,M_1\mid M_2\,]$, and $[\,M_1^\top\mid M_2^\top\,]$. Hence $(M_1',M_2')\in\variety{I'_{\mathcal A}}$. Since $n_1\leq N$, Theorem~\ref{thm:set-theoretic-small} gives $(M_1',M_2')\in\mathcal U'_{\mathcal A}$.

Let $P\in\CC^{n_2\times N}$ select the first $n_2$ rows and let $Q\in\CC^{N\times n_0}$ select the first $n_0$ columns:
\[P = \begin{bmatrix}
    I_{n_2} \mid  0_{n_2\times(N-n_2)}
\end{bmatrix},\qquad Q = \begin{bmatrix}
    I_{n_0} \\ 0_{(N-n_0)\times n_0}
\end{bmatrix}.\]
Define $\pi(X,Y)=(PXQ,PYQ)$. Then $\pi(M_1',M_2')=(M_1,M_2)$. We claim that $\pi(\im\varphi'_{\mathcal A})=\im\varphi_{\mathcal A}$. Indeed, for $i=1,2$,
$$
P\bigl(W^{(2)}\diag(A_i)W^{(1)}\bigr)Q
=
(PW^{(2)})\diag(A_i)(W^{(1)}Q),
$$
so $\pi(\im\varphi'_{\mathcal A})\subseteq\im\varphi_{\mathcal A}$. Conversely, any $W^{(2)}\in\CC^{n_2\times n_1}$ and $W^{(1)}\in\CC^{n_1\times n_0}$ can be extended by zero rows and columns to matrices of sizes $N\times n_1$ and $n_1\times N$, respectively, giving the reverse inclusion. Since $\pi$ is linear, it is Zariski continuous. Therefore
$$
(M_1,M_2)\in
\pi(\mathcal U'_{\mathcal A})
\subseteq
\overline{\pi(\im\varphi'_{\mathcal A})}
=
\overline{\im\varphi_{\mathcal A}}
=
\mathcal U_{\mathcal A}.
$$
Thus $\variety{I_{\mathcal A}}\subseteq\mathcal U_{\mathcal A}$, completing the proof.
\end{proof}

\section{Defining ideals}\label{sec:ideal-theoretic}

In this section, we build off of the set-theoretic results of the previous section to show that the ideal $I_{\mathcal{A}}$ is radical and is therefore the defining ideal of the two-block ReLU pattern variety. 
By Theorem~\ref{thm:set-theoretic-general}, we know that $I_{\mathcal{A}}$ set-theoretically defines an irreducible variety, and so the results of this section show that it is a prime ideal.

\subsection{Square case}\label{subsec:Isquare}

We first consider the square case $n_0=n_1=n_2=n$. Recall that
$n=a+b+c$. In this case,
$$I_{\mathcal A}=I_{n-b+1}(M_1)+I_{n-a+1}(M_2)+I_{n-c+1}(M_1-M_2).$$
We will show that $I_{\mathcal A}$ is prime in this case by using Serre's criterion. In particular, we will show that $R/I_{\mathcal{A}}$ is Cohen--Macaulay and generically reduced; these together imply that $I_{\mathcal{A}}$ is radical. 

We begin by proving that $R/I_{\mathcal{A}}$ is Cohen--Macaulay.

\begin{lemma}\label{lem:square-CM}
The ring $R/I_{\mathcal A}$ is Cohen--Macaulay and
$$
\dim R/I_{\mathcal A}
=
2n^2-a^2-b^2-c^2.
$$
\end{lemma}
\begin{proof}
Let $Z=(z_{ij})$ be a generic $n\times n$ matrix and set
$$S=\CC[M_1,M_2,Z],\qquad
K=
I_{n-b+1}(M_1)
+
I_{n-a+1}(M_2)
+
I_{n-c+1}(Z).$$
Since the three determinantal ideals involve disjoint sets of variables, we have
$$
S/K
\cong
\frac{\CC[M_1]}{I_{n-b+1}(M_1)}
\otimes_{\CC}
\frac{\CC[M_2]}{I_{n-a+1}(M_2)}
\otimes_{\CC}
\frac{\CC[Z]}{I_{n-c+1}(Z)}.
$$
Each factor is a generic determinantal ring and hence Cohen--Macaulay. Therefore $S/K$ is Cohen-Macaulay. The three determinantal varieties have codimensions $b^2$, $a^2$, and $c^2$, respectively, so $\dim S/K=3n^2-a^2-b^2-c^2.$

Let $L\subseteq S/K$ be the ideal generated by the $n^2$ linear forms
$z_{ij}-x_{ij}+y_{ij}$. Then $(S/K)/L\cong R/I_{\mathcal A}.$ By Theorem~\ref{thm:set-theoretic-general},
$\variety{I_{\mathcal A}}=\mathcal U_{\mathcal A}$, and the dimension
formula for $\mathcal U_{\mathcal A}$ from~\cite{alexandr_montufar_2027}
gives
\[
\dim R/I_{\mathcal A}
=
2n^2-a^2-b^2-c^2.
\]
Thus quotienting $S/K$ by the $n^2$ generators of $L$ lowers the dimension
by exactly $n^2$. Since $S/K$ is Cohen--Macaulay, an ideal of codimension $n^2$ and generated by
$n^2$ elements is generated by a regular
sequence. Hence the generators of $L$ form a regular sequence. Quotienting
a Cohen--Macaulay ring by $L$ preserves the Cohen--Macaulay
property, so $(S/K)/L\cong R/I_{\mathcal A}$ is Cohen--Macaulay, as~desired.
\end{proof}

Recall that $D_1=\diag(I_a,0_b,I_c)$ and $D_2=\diag(0_a,I_b,I_c)$. The point $(D_1,D_2)$ is in the variety $\mathcal{U}_{\mathcal{A}}$, and so it can be considered as a point in $\Spec(R/I_{\mathcal{A}})$. We show that $R/I_{\mathcal{A}}$ is generically reduced by showing that $(D_1,D_2)$ is a regular point.

\begin{lemma}\label{lem:square-generically-reduced}
At the point $(D_1,D_2)$, the Zariski tangent space to
$\Spec(R/I_{\mathcal A})$ has dimension
$$
\dim T_{(D_1,D_2)}\Spec(R/I_{\mathcal A})
=
2n^2-a^2-b^2-c^2.
$$
Consequently, $R/I_{\mathcal A}$ is generically reduced.
\end{lemma}

\begin{proof}
Let $\Sigma_r\subseteq\CC^{n\times n}$ be the determinantal
variety defined by the $(r+1)$-minors. For a matrix $M$ of rank
exactly $r$, its tangent space is
\[
T_M\Sigma_r
=
\{H\in\CC^{n\times n}:H(\ker M)\subseteq\im M\};
\]
see~\cite[Example~14.16]{harris}. Write $\CC^n=\CC^a\oplus\CC^b\oplus\CC^c$; for a pair of matrices $(H_1,H_2)$, express each~$H_i$ in the corresponding block form
$H_i=((H_i)_{jk})_{j,k=1}^3$. The matrices
\[
D_1=\diag(I_a,0_b,I_c),
\qquad
D_2=\diag(0_a,I_b,I_c),
\qquad
D_1-D_2=\diag(I_a,-I_b,0_c),
\]
have ranks $n-b$, $n-a$, and $n-c$, respectively.
Since
\[
I_{\mathcal A}
=
I_{n-b+1}(M_1)
+
I_{n-a+1}(M_2)
+
I_{n-c+1}(M_1-M_2),
\]
the tangent space to $\Spec(R/I_{\mathcal A})$ at $(D_1,D_2)$
is obtained by imposing the three determinantal tangent space
conditions simultaneously. Applying the formula above to
$D_1$, $D_2$, and $D_1-D_2$, with tangent directions
$H_1$, $H_2$, and $H_1-H_2$, respectively, gives
\[
T_{(D_1,D_2)}\Spec(R/I_{\mathcal A})
=
\bigl\{(H_1,H_2):
(H_1)_{22}=0,\quad
(H_2)_{11}=0,\quad
(H_1)_{33}=(H_2)_{33}
\bigr\}.
\]
These conditions impose $b^2$, $a^2$, and $c^2$ independent
linear equations, respectively, since the three groups involve
disjoint sets of coordinates. Therefore
\[
\dim T_{(D_1,D_2)}\Spec(R/I_{\mathcal A})
=
2n^2-a^2-b^2-c^2.
\]

Set $X=\Spec(R/I_{\mathcal A})$ and $p=(D_1,D_2)$.
By Theorem~\ref{thm:set-theoretic-general},
$\sqrt{I_{\mathcal A}}=J_{\mathcal A}$.
Since $J_{\mathcal A}$ is prime, $X$ is irreducible.
Moreover, Lemma~\ref{lem:square-CM} gives $\dim X=2n^2-a^2-b^2-c^2.$ As $X$ is an irreducible scheme of finite type over $\CC$
and $p$ is a closed point, the local dimension of $X$ at $p$
equals $\dim X$. The tangent space computation therefore
shows that $p$ is a regular point of $X$.
The regular locus of $X$ is open and contains $p$, so it is
nonempty and hence dense by irreducibility.
Thus $X$ is reduced on a dense open subset, and
$R/I_{\mathcal A}$ is generically reduced.
\end{proof}

We are now ready to prove primeness in the square case.

\begin{theorem}\label{thm:square-prime}
    $I_\cA$ is prime when $n_0=n_1=n_2$.
\end{theorem}

\begin{proof}
By Lemma~\ref{lem:square-CM}, the ring $R/I_{\mathcal A}$
is Cohen--Macaulay and therefore has no embedded associated
primes. By Lemma~\ref{lem:square-generically-reduced}, it is
generically reduced. Hence $R/I_{\mathcal A}$ is reduced by Serre's criterion $(R_0)+(S_1)$,
so $I_{\mathcal A}=\sqrt{I_{\mathcal A}}$.
Theorem~\ref{thm:set-theoretic-general} gives
$\sqrt{I_{\mathcal A}}=J_{\mathcal A}$.
Since $J_{\mathcal A}$ is prime, so is $I_{\mathcal A}$.
\end{proof}

\begin{remark}
    We note that Lemma~\ref{lem:square-CM} in fact shows that $R/I_{\mathcal A}$ is Gorenstein in the square case. Similar arguments to those in this section can be used to show that $R/I_{\mathcal{A}}$ is always Gorenstein (and therefore Cohen--Macaulay) whenever $n_0=n_2\leq n_1$.
    Computations in \texttt{Macaulay2} suggest that the ring is always Cohen--Macaulay for arbitrary choice of layer dimensions $n_0,n_1,n_2$ \cite{M2}. 
\end{remark}

\subsection{Reduction for smaller input and output dimensions}

The following lemma allows us to deduce from Theorem~\ref{thm:square-prime} that $I_{\mathcal{A}}$ is prime when the input and output dimensions of the architecture are no more than the dimension of the hidden layer.

\begin{lemma}\label{lem:ideal red to smaller arch}
    Let $\cA,\cB$ be two pairs of activation patterns for networks with architectures $(n_0,n_1,n_2)$ and $(m_0,n_1,m_2)$ respectively, where $m_0\leq n_0$ and $m_2\leq n_2$, and the values of $a,b,c$ are the same for $\cA$ and $\cB$. 
    If $I_{\cA}$ is prime, then so is $I_{\cB}$.
\end{lemma}

\begin{proof}
    Let $\pi:\CC^{n_2\times n_0}\to\CC^{m_2\times m_0}$ be the linear map taking the upper-left $(m_2\times m_0$)-submatrix of an $(n_2\times n_0)$-matrix, and let $\iota:\CC^{m_2\times m_0}\to\CC^{n_2\times n_0}$ be the linear map taking an $(m_2\times m_0)$-matrix and adding $0$ entries to obtain an $(n_2\times n_0)$-matrix. Note that $\pi\circ\iota=\id$. 

    Let $R_\cA$ and $R_\cB$ denote the respective polynomial rings. We have induced pullback maps on these rings: $\pi^*:R_\cB\to R_\cA$ and $\iota^*:R_\cA\to R_\cB$. The former map is inclusion of rings, so $\pi^*(I_\cB)\subseteq I_\cA$. The latter map sends variables to $0$. Every defining minor of $I_\cA$ is then mapped to either a defining minor of $I_\cB$ or to $0$, so $\iota^*(I_\cA)=I_\cB$.

    We claim that $I_\cB=(\pi^*)^{-1}(I_\cA)$. The containment $I_\cB\subseteq (\pi^*)^{-1}(I_\cA)$ follows from $\pi^*(I_\cB)\subseteq I_\cA$. For the reverse containment, assume that $\pi^*(f)\in I_\cA$. Then
    \[f=\id^*(f)=(\iota^*\circ\pi^*)(f)\in\iota^*(I_\cA)=I_\cB.\]
    Therefore, $I_\cB$ is the preimage under a ring homomorphism of $I_\cA$. If $I_\cA$ is prime, then $I_\cB$ is prime.
\end{proof}

\subsection{Inheritance to larger input or output dimensions}\label{subsec:inheritance}

To show that $I_{\mathcal{A}}$ is prime for general architectures $(n_0,n_1,n_2)$, it remains to consider the case where $n_0,n_2\geq n_1$. 
By Lemma~\ref{lem:ideal red to smaller arch}, it suffices to
prove primality when both the input and output dimensions are
at least $n_1$. Indeed, an arbitrary architecture $(n_0,n_1,n_2)$
is obtained by restriction from the
$(\max\{n_0,n_1\},n_1,\max\{n_2,n_1\})$ architecture.

Let $\cB$ be a pair of activation patterns with such an architecture and values $a,b,c$; let $\cA$ be a pair of activation patterns for the $(n_1,n_1,n_1)$ architecture, same $a,b,c$ values. We already know $I_\cA$ is prime, and we want to show that $I_\cB$ is prime.

Let $\iota:\CC^{n_1\times n_1}\to \CC^{n_2\times n_0}$ be the inclusion map, and let $\iota^*:R_\cB\to R_\cA$ be the induced map on coordinate rings. 

\begin{lemma}\label{lem:ideal pullback equality}
   $I_{\cA}=\iota^*(I_{\cB})=\iota^*(\sqrt{I_\cB})$.
\end{lemma}

\begin{proof}
    As in the proof of Lemma~\ref{lem:ideal red to smaller arch}, we have that $\iota^*(I_\cB)=I_\cA$.

    Suppose $f\in\sqrt{I_\cB}$, so $f^m\in I_\cB$ for some $m$. Then $\iota^*(f^m)=\iota^*(f)^m\in I_\cA$, by the first equality. We have already shown that $I_\cA$ is prime in \S\ref{subsec:Isquare}, as $\cA$ is the square matrix architecture. Therefore, $\iota^*(f)\in I_\cA$ and $\iota^*(\sqrt{I_\cB})\subseteq I_\cA$. Since $I_\cB\subseteq\sqrt{I_\cB}$, the reverse containment also holds.
\end{proof}

Let $E,F$ be $\CC$-vector spaces with $E\cong\CC^{n_2}$ and $F\cong\CC^{n_0}$ so that the coordinate ring $R_{\cB}$ of pairs of $(n_2\times n_0)$-matrices can be identified with the symmetric algebra $\Sym(E\otimes F\otimes \CC^2)$. The group $\GL(E)\times\GL(F)$ naturally acts on $R_\cB$. The following lemma decomposes $R_\cB$ and ideals related to $I_\cB$ into Schur modules, which are the simple polynomial representations with respect to the group $\GL(E)\times \GL(F)$. For further background on Schur decompositions, we refer the reader to \cite{BCRV,landsberg}.

\begin{lemma}\label{lem:schur decomp}
    Let $I'_\cB$ denote the ideal generated by $(n_1+1)$-minors of $[X\mid Y]$ and $[X^\top \mid Y^\top]$. Then for any $d\geq 0$, we have the following decompositions as $\GL(E)\times\GL(F)$-representations of the spaces of degree $d$ polynomials of the ring and ideal:
    \begin{align*}
        (R_\cB)_d&=\bigoplus_{\lambda,\mu\vdash d}\bbS_\lambda(E)\otimes \bbS_\mu(F)\otimes M_{\lambda,\mu},\\
        (I'_\cB)_d&=\bigoplus_{\substack{\lambda,\mu\vdash d,\\\ell(\lambda)>n_1\\\text{ or $\ell(\mu)>n_1$}}}\bbS_\lambda(E)\otimes \bbS_\mu(F)\otimes M_{\lambda,\mu},
    \end{align*}
    where $M_{\lambda,\mu}$ is the same in both lines and is independent of $n_0,n_2$. Furthermore, these decompositions are functorial in $E$ and $F$.
\end{lemma}

\begin{proof}
    For the decomposition of $(R_\cB)_d$, we first identify this space with $\Sym^d(E\otimes(F\otimes \CC^2))$, and then apply Cauchy's identity:
    \[(R_\cB)_d=\Sym^d(E\otimes (F\otimes \CC^2))=\bigoplus_{\lambda\vdash d}\bbS_\lambda(E)\otimes \bbS_\lambda(F\otimes \CC^2).\]
    As a $\GL(F)$-module, $\bbS_\lambda(F\otimes \CC^2)$ can be further decomposed as
    \[\bbS_\lambda(F\otimes \CC^2)=\bigoplus_{\mu\vdash d}\bbS_\mu(F)\otimes M_{\lambda,\mu},\]
    where the dimension of $M_{\lambda,\mu}$ as a vector space depends only on $\lambda,\mu$ via Kronecker coefficients.

    For the decomposition of $I'_\cB$, first note that this ideal is stable under the action of $\GL(E)\times\GL(F)$. Identify the concatenation matrix $[X\mid Y]$ with $E\otimes (F\otimes \CC^2)$. Then the ideal of $(n_1+1)$-minors of $[X\mid Y]$ in degree $d$ has the Schur decomposition
    \[\bigoplus_{\substack{\lambda\vdash d,\\ \ell(\lambda)>n_1}}\bbS_\lambda(E)\otimes \bbS_\lambda(F\otimes \CC^2)=\bigoplus_{\substack{\lambda,\mu\vdash d,\\ \ell(\lambda)>n_1}}\bbS_\lambda(E)\otimes\bbS_\mu(F)\otimes M_{\lambda,\mu}.\]
    Similarly, the matrix $[X^\top\mid Y^\top]$ can be identified with $F\otimes (E\otimes \CC^2)$, so the ideal of its $(n_1+1)$-minors in degree $d$ has the Schur decomposition
    \[\bigoplus_{\substack{\lambda,\mu\vdash d,\\ \ell(\mu)>n_1}}\bbS_\lambda(E)\otimes\bbS_\mu(F)\otimes M_{\lambda,\mu}.\]
    Together, these decompositions prove the description of $(I'_\cB)_d$.
\end{proof}

\begin{theorem}
    The ideal $I_\cB$ is prime.
\end{theorem}

\begin{proof}
    Let $P_\cB=\sqrt{I_\cB}$ denote the radical of $I_\cB$. We already know from Theorem~\ref{thm:set-theoretic-general} that $P_\cB$ is prime. These ideals are both homogeneous, so it suffices to show that $(P_\cB)_d=(I_\cB)_d$ are equal in each degree $d\geq 0$.

    Both $I_\cB$ and $P_\cB$ are $\GL(E)\times \GL(F)$-stable, so their degree $d$ pieces decompose as $\GL(E)\times\GL(F)$-representations:
    \[(I_\cB)_d=\bigoplus_{\lambda,\mu\vdash d}\bbS_\lambda(E)\otimes\bbS_\mu(F)\otimes N^I_{\lambda,\mu},\qquad (P_\cB)_d=\bigoplus_{\lambda,\mu\vdash d}\bbS_\lambda(E)\otimes\bbS_\mu(F)\otimes N^P_{\lambda,\mu},\]
    where $N^I_{\lambda,\mu}\subseteq N^P_{\lambda,\mu}\subseteq M_{\lambda\mu}$ for all pairs of partitions $\lambda,\mu$; this follows from the fact that each $\bbS_\lambda(E)\otimes\bbS_\mu(F)$ is a simple $\GL(E)\times\GL(F)$-representation. We want to show that the multiplicity spaces $N^I_{\lambda,\mu}=N^P_{\lambda,\mu}$ are equal for all pairs of partitions.

    The pullback $\iota^*$ acts on $E\otimes F\otimes \CC^2$ by $\pi_E\otimes \pi_F\otimes \id$, where $\pi_E:E\to E_0$ and $\pi_F:F\to F_0$ are projections onto $n_1$-dimensional subspaces $E_0$ and $F_0$. Therefore, $\iota^*$ acts on each summand of a Schur decomposition by replacing $E,F$ with $E_0,F_0$, and leaving the multiplicity space unchanged. Since $E_0,F_0$ are $n_1$-dimensional, we have that $\bbS_\lambda(E_0)\otimes\bbS_\mu(F_0)=0$ if either $\ell(\lambda)>n_1$ or $\ell(\mu)>n_1$. Therefore, we have Schur decompositions
    \[(\iota^*I_\cB)_d=\bigoplus_{\substack{\lambda,\mu\vdash d\\\ell(\lambda),\ell(\mu)\leq n_1}}\bbS_\lambda(E_0)\otimes\bbS_\mu(F_0)\otimes N^I_{\lambda,\mu},\qquad (\iota^*P_\cB)_d=\bigoplus_{\substack{\lambda,\mu\vdash d\\\ell(\lambda),\ell(\mu)\leq n_1}}\bbS_\lambda(E_0)\otimes\bbS_\mu(F_0)\otimes N^P_{\lambda,\mu}.\]

    By Lemma~\ref{lem:ideal pullback equality}, we have that
    \[(\iota^*I_\cB)_d=(\iota^*P_\cB)_d=(I_\cA)_d,\]
    so $N^I_{\lambda,\mu}=N^P_{\lambda,\mu}$ for all pairs of partitions with $\ell(\lambda),\ell(\mu)\leq n_1$.

    If $\ell(\lambda)>n_1$ or $\ell(\mu)>n_1$, then Lemma~\ref{lem:schur decomp} implies that $N^I_{\lambda,\mu}=M_{\lambda,\mu}$, forcing $N^I_{\lambda,\mu}=N^P_{\lambda,\mu}$ as well. Since $N^I_{\lambda,\mu}=N^P_{\lambda,\mu}$ for all pairs of partitions, we see that $I_\cB=P_\cB$ is prime, as desired.
\end{proof}

Applying this theorem to the enlarged architecture
$(\max\{n_0,n_1\},n_1,\max\{n_2,n_1\})$ and then using
Lemma~\ref{lem:ideal red to smaller arch} proves primality of $I_{\mathcal A}$
for arbitrary input and output dimensions.
Together with Theorem~\ref{thm:set-theoretic-general},
this gives $I_{\mathcal A}=J_{\mathcal A}$ and completes
the proof of Theorem~\ref{thm:main}.

\section{An incidence variety}\label{sec:incidence}

In this section, we give an alternative proof of irreducibility of the variety defined by $I=I_\cA$ in the square case of $n_0=n_1=n_2$. This section is independent of the results in Section~\ref{sec:matrix pencils}, and uses an incidence variety to analyze $\cV(I)$. 
By continuing with the arguments of Section~\ref{sec:ideal-theoretic}, one can again conclude that $I_\cA$ is prime for arbitrary architectures. 

Set $n=n_0=n_1=n_2$, and let $V=\mathbb{C}^n$. Identify $\mathrm{End}(V)\times\mathrm{End}(V)\cong\mathbb{C}^{2n^2}$, and let $\mathcal{G}=\mathrm{Gr}(b,V)\times\mathrm{Gr}(a,V)\times\mathrm{Gr}(c,V)$. We define an incidence variety $\mathcal{Z}\subseteq \mathbb{C}^{2n^2}\times\mathcal{G}$ as the Zariski closure of the following set:
    \[\{(X,Y,A,B,C):XA=0,YB=0,(X-Y)C=0\}.\]
    Let $\pi_1$ and $\pi_2$ denote the projections of $\mathcal{Z}$ to $\mathbb{C}^{2n^2}$ and to $\mathcal{G}$ respectively. By construction, $\im\pi_1=\mathcal{V}(I)$; here, we use the fact that $a+b+c=n$, so the nullity of $X$ is at least $n-(a+c)=b$. 
    
    Let $U\subseteq \mathcal{G}$ denote the subset
    \[U=\{(A,B,C):A\oplus B\oplus C = V\}\subseteq\mathcal{G}.\]

\begin{lemma}\label{lem:open subset U dim} The following statements hold:
    \begin{enumerate}[label={(\arabic*)}]
        \item $U\subseteq\mathcal{G}$ is a dense, open, irreducible subset.
        \item $\pi_2^{-1}(U)\subseteq\mathcal{Z}$ is a vector bundle of rank $n^2$ over $U$; it is irreducible of dimension $2n^2-a^2-b^2-c^2$.
    \end{enumerate}
\end{lemma}

\begin{proof}
We prove the two statements. To show (1): since $\mathcal{G}$ is irreducible, it suffices to show that $U$ is open and nonempty. Nonemptiness follows from picking a basis of $V$ and partitioning it to get an $A,B,C$. To show openness: pick a basis of $V$, write $A,B,C$ with respect to this basis, make an $n\times n$ matrix with these coefficients. Then $(A,B,C)\in U$ if and only if this matrix is invertible, which is an open condition.

For (2), let $(A,B,C)\in U$. We show that $\pi_2^{-1}(A,B,C)$ is isomorphic as vector spaces to
        \[\mathrm{Hom}_{\mathbb{C}}(A,V)\oplus\mathrm{Hom}_{\mathbb{C}}(B,V)\oplus\mathrm{Hom}_{\mathbb{C}}(C,V),\]
        which has dimension
        \[an+bn+cn=n^2.\]
        Given $(X,Y)\in\pi_2^{-1}(A,B,C)$, write $X,Y$ with respect to a basis given by $A\oplus B\oplus C=V$. Then $(X-Y)C=0$ means that $X,Y$ agree on $C$, determining a linear map $C\to V$. A linear map $B\to V$ is determined by $X$, and a linear map $A\to V$ is determined by $Y$. This correspondence is reversible since $XA=0$ and $YB=0$. The fibers vary algebraically with $(A,B,C)$. This shows that $\pi_2^{-1}(U)$ is a vector bundle of rank $n^2$, and so its dimension is
        \begin{align*}
            \dim\pi_2^{-1}(U)&=n^2+\dim\mathcal{G}\\
            &=n^2+\dim\mathrm{Gr}(a,V)+\dim\mathrm{Gr}(b,V)+\dim\mathrm{Gr}(c,V)\\
            &=n^2+a(n-a)+b(n-b)+c(n-c)\\
            &=2n^2-a^2-b^2-c^2.
        \end{align*}
        Finally, since $U$ is irreducible, so is $\pi_2^{-1}(U)$.
\end{proof}

Let $\partial \mathcal{G}=\mathcal{G}\setminus U$. The following lemma bounds the dimension of its preimage under $\pi_2$; we will then use this to show that $\pi_1(\pi_2^{-1}(U))$ is dense in $\mathcal{V}(I)$.

\begin{lemma}\label{lem:pi_2 partial G dim}
    We have $\dim\pi_2^{-1}(\partial\mathcal{G})<2n^2-a^2-b^2-c^2$.
\end{lemma}

\begin{proof}
    Fix $(A,B,C)\in\mathcal{G}$, and let $t=\dim(A\cap B\cap C)$. We show that the dimension of the fiber over this point is $\dim\pi_2^{-1}(A,B,C)=n^2+nt$. Consider the following linear map
        \begin{align*}
            \varphi_{A,B,C}:\mathrm{End}(V)^2&\to\mathrm{Hom}_\CC(A,V)\oplus\mathrm{Hom}_\CC(B,V)\oplus\mathrm{Hom}_\CC(C,V)\\
            (X,Y)&\mapsto (X|_A, Y|_B, (X-Y)|_C).
        \end{align*}
        Then $\pi_2^{-1}(A,B,C)=\ker\varphi_{A,B,C}$. To compute the rank of $\varphi_{A,B,C}$, we use the following linear map
        \[\psi:V^*\oplus V^*\to A^*\oplus B^*\oplus C^*,\qquad (f,g)\mapsto (f|_A,g|_B,(f-g)|_C),\]
        and its dual:
        \[\psi^*:A\oplus B\oplus C\to V\oplus V,\qquad (a,b,c)\mapsto (a+c, b-c).\]
        Then $\ker\psi^*=\{(-v,v,v):v\in A\cap B\cap C\}$, so $\psi$ has rank $n-t$. By identifying $\mathrm{End}(V)=V^*\otimes V$ and similarly for $A^*,B^*,C^*$, the original map $\varphi_{A,B,C}$ can be obtained from $\psi$ by tensoring with the identity map $\mathrm{id}:V\to V$, and so the rank of $\varphi_{A,B,C}$ is $n(n-t)$. Therefore,
        \[\dim\pi_2^{-1}(A,B,C)=2n^2-n(n-t)=n^2+nt.\]

        For $t\geq 1$, let $\mathcal{G}_{\geq t}=\{(A,B,C):\dim(A\cap B\cap C)\geq t\}$. A point in $\mathcal{G}_{\geq t}$ can be given by choosing a $t$-subspace of $V$, and then choosing appropriate $A,B,C$ containing the $t$-subspace. Therefore,
        \begin{align*}
            \dim\mathcal{G}_{\geq t}&\leq t(n-t)+(a-t)(n-a)+(b-t)(n-b)+(c-t)(n-c)\\
            &=\dim\mathcal{G}-t(n+t).
        \end{align*}
        Combining this with the dimension of $\pi_2^{-1}(A,B,C)$, we see that the inverse image under $\pi_2$ of $\mathcal{G}_{\geq t}\setminus{\mathcal{G}_{\geq t+1}}$ has dimension at most
        \[\dim\mathcal{G}-t(n+t)+n^2+nt=2n^2-a^2-b^2-c^2-t^2<2n^2-a^2-b^2-c^2.\]

        It remains to consider the subset of $\partial \mathcal{G}$ given by
        \[\mathcal{G}_0=\{(A,B,C)\in\partial\mathcal{G}:\dim(A\cap B\cap C)=0\}.\]
        The dimension of $\pi_2^{-1}(A,B,C)$ for $(A,B,C)\in\mathcal{G}_0$ is $n^2+n(0)=n^2$. Furthermore, $\dim\partial\mathcal{G}<\dim\mathcal{G}$, and so
        \[\dim\pi_2^{-1}(\mathcal{G}_0)<n^2+\dim\mathcal{G}=2n^2-a^2-b^2-c^2.\]
\end{proof}

\begin{lemma}\label{lem:dense open subset V(I) dim}
    $\pi_1(\pi_2^{-1}(U))$ is dense in $\mathcal{V}(I)$, and $\mathcal{V}(I)$ is irreducible. Furthermore, $\dim\mathcal{V}(I)=2n^2-a^2-b^2-c^2$. 
\end{lemma}

\begin{proof}
     Let $X=\pi_1(\pi_2^{-1}(\partial\mathcal{G}))$ and $Y=\pi_1(\pi_2^{-1}(U))$. Then $X\cup Y=\im\pi_1=\mathcal{V}(I)$. Let $Z\subseteq\mathcal{V}(I)$ be an irreducible component. In order to show that $Y$ is dense in $\mathcal{V}(I)$, it suffices to show that $Z=\overline{Y}$.

        As in the proof of Lemma~\ref{lem:square-CM}, we have that $R/I$ is isomorphic to a Cohen--Macaulay ring of dimension $3n^2-a^2-b^2-c^2$, modulo an ideal generated by $n^2$ linear forms. By Krull height theorem, this implies that every irreducible component has dimension at least $2n^2-a^2-b^2-c^2$. But by Lemma~\ref{lem:pi_2 partial G dim}, $\dim X<2n^2-a^2-b^2-c^2$, so $Z\not\subseteq X$. In particular, $Z\setminus X\subseteq Y$ is a nonempty subset of $Y$. Properness of $\pi_1$ means that $X$ is closed, so $Z\setminus X$ is a dense open subset in $Z$. Altogether, this implies that $Z\subseteq\overline{Y}$. From Lemma~\ref{lem:open subset U dim}, we see that $\dim\overline{Y}\leq 2n^2-a^2-b^2-c^2$, and so we must have that $Z=\overline{Y}$.

        Since $Z$ is an arbitrary irreducible component, this all implies that $\mathcal{V}(I)=\overline{Y}$, so $Y=\pi_1(\pi_2^{-1}(U))$ is dense in $\mathcal{V}(I)$. Therefore, $\mathcal{V}(I)$ is irreducible and has dimension $2n^2-a^2-b^2-c^2$.
\end{proof}

\begin{theorem}
    $I$ is prime.
\end{theorem}

\begin{proof}
    By Lemma~\ref{lem:dense open subset V(I) dim}, $\dim\mathcal{V}(I)=2n^2-a^2-b^2-c^2$ and $\mathcal{V}(I)$ is irreducible. Using the same argument as in Lemmas~\ref{lem:square-CM} and \ref{lem:square-generically-reduced}, $R/I$ is Cohen--Macaulay and generically reduced. Altogether, this implies that $I$ is prime.
\end{proof}

\section{Conclusion and future work}\label{sec:conclusion}

\subsection{Towards defining ideals for more activation regions}\label{subsec:more act regions}

For a single activation pattern, the pattern variety is a classical determinantal variety. Our results give a complete description for two patterns: the minors expressing the rank bounds on the two blocks, their difference, and their concatenations generate the defining ideal. 

When the selected activation patterns include those of \textit{all} full-dimensional activation regions for a fixed choice of parameters, the defining ideal admits another explicit description. For pairwise distinct activation hyperplanes, this description is established in~\cite{finkel2026canonicalactivationpatternsrelu}, as part of a broader study of the relation between the combinatorics of activation patterns and the geometry of their pattern varieties.
Fix a collection $\mathcal A$ containing the patterns of all full-dimensional regions of an arrangement of $n_1$ pairwise distinct activation hyperplanes, and let $P$ be the matrix with these patterns as rows.
For each neuron $i$, there are two adjacent regions whose patterns differ by $e_i$. 
Hence the row span of $P$ contains every standard basis vector, so $\rank P=n_1$. 
Now consider the linear map
$$
\Phi_{\mathcal A}:(N_1,\ldots,N_{n_1})
\longmapsto
\left(\sum_{i=1}^{n_1}A_iN_i\right)_{A\in\mathcal A},
\qquad N_i\in\CC^{n_2\times n_0}.
$$
Since $P$ has full column rank, $\Phi_{\mathcal A}$ is a linear embedding. 
The network parametrization has image $\Phi_{\mathcal A}(D_1^{n_1})$, where $D_1$ is the variety of matrices of rank at most one: each $N_i$ is the outer product of column $i$ of $W^{(2)}$ and row $i$ of $W^{(1)}$. 
Since $\Phi_{\mathcal A}$ is a linear embedding, this image is closed. Thus~$\Phi_{\mathcal A}$ restricts to an isomorphism between $D_1^{n_1}$ and the pattern variety. In block coordinates, the defining ideal is generated by the linear equations of $\operatorname{im}\Phi_{\mathcal A}$ and the $2\times2$ minors of the matrices $N_i$, each expressed as a difference of two adjacent blocks. 
More generally, \cite{finkel2026canonicalactivationpatternsrelu}
covers every collection that contains, for each neuron, two patterns differing
only in that neuron, and in which any two patterns are connected by a sequence of
such single-neuron changes.
The argument above uses only that $P$ has full column rank, so it applies to
every collection with this property. The matrices $N_i$ are recovered as linear combinations of the blocks, not necessarily as differences of adjacent blocks.

The remaining challenge concerns proper subcollections of three or more activation patterns whose neuron contributions cannot be recovered linearly, even after grouping neurons that are active in exactly the same selected patterns.
In these cases, eliminating the shared rank-one summands can introduce higher-degree compatibility equations involving several blocks at once, beyond those imposed by individual pairs.
The quintic generators in~\cite[Example~37]{alexandr_montufar_2027} illustrate this phenomenon. 
Whether such determinantal constructions generate the defining ideal is also not known; see~\cite[Example~38]{alexandr_montufar_2027}. 
It remains to determine which collections of patterns admit a determinantal description and how their defining equations depend on the overlaps among the active neurons.

\subsection{Gr\"obner bases}

A natural question to ask about the defining ideal $I_\cA$ is whether or not its natural generating set forms a Gr\"obner basis with respect to some term order. The following proposition demonstrates that this is not generally the case.

\begin{proposition}\label{prop:gb countereg}
    Suppose $n_0=n_2=5$, $n_1=3$, and $a=b=c=1$. Then the natural generating set of $I_\cA$ given in Theorem~\ref{thm:main} does not form a Gr\"obner basis with respect to any term order.
\end{proposition}

\begin{proof}
    Fix a term order $\prec$. We first show that there exists a $3$-minor of $X-Y$, indexed by row set $S\subseteq [5]$ and column set $T\subseteq[5]$, for which its leading term consists of only $x$ variables or of only $y$ variables. Let $D=\det(X-Y)$, and let $m=\LT(D)$ be the leading term of this determinant. This monomial $m$ consists of $5$ distinct variables; by pigeonhole principle, at least $3$ of them are all $x$ variables, or all $y$ variables. Choose $3$ such variables, and denote their product as $m_0$. Let $S,T\subseteq[5]$ denote the row and column indices respectively of the variables appearing in $m_0$, so $\#S=\#T=3$. We claim that $m_0=\LT(\det(X-Y)_{S,T})$.

    Suppose a larger monomial $s\succ m_0$ is a term in $\det(X-Y)_{S,T}$. Then $s(m/m_0)$ would be a term appearing in $D$ with nonzero coefficient, and $s(m/m_0)\succ m=m_0(m/m_0)$ by multiplicativity of the term order. This contradicts the assumption that $m=\LT(D)$.

    For the row and column sets $S$ and $T$ just found, consider the following three polynomials which are in the natural generating set of $I_\cA$:
    \[f=\det X_{S,T},\quad g = \det Y_{S,T},\quad h = \det (X-Y)_{S,T}.\]
    Consider the polynomial $p=h-f+g\in I_\cA$. We show that the leading term $\LT(p)$ of $p$ is not generated by the leading terms of the generating set. First note $\deg p = 3$, that every monomial of $p$ involves both $x$ and $y$ variables, and that $p$ is nonzero. Therefore, $\LT(p)$ cannot be divisible by any leading term of a $3$-minor of $X$ or of $Y$, or of a $4$-minor of either of the concatenation matrices. It remains to show that $\LT(p)$ is not the leading term of any $3$-minor of $X-Y$. Given the choice of row and column indices, if $\LT(p)$ were to be the leading term of a $3$-minor of $X-Y$, then it must be equal to $\LT(h)=\LT(\det(X-Y)_{S,T})$. However, we assumed that the leading term of $h$ is not divisible by both $x$ and $y$ variables. This concludes the proof.
\end{proof}

This raises the problem of finding a Gr\"obner basis of $I_\cA$ with respect to a natural term order for any choice of activation pattern $\cA$ and network architecture.

\bibliographystyle{plain}
\bibliography{references}

\authorinfo

\end{document}